\documentclass[11pt]{article}
\usepackage{amsmath, amssymb, amscd}
\usepackage[dvips]{graphics}

\newtheorem{thm}{Theorem}
\newtheorem{prop}[thm]{Proposition}
\newtheorem{lem}[thm]{Lemma}

\newtheorem{rem}[thm]{Remark}
\newtheorem{df}[thm]{Definition}

\renewcommand{\epsilon}{\varepsilon}
\renewcommand{\phi}{\varphi}

\renewcommand{\P}{\operatorname{P}}

\newcommand{\BB}{\mathbb}

\newcommand{\pf}{\noindent {\it Proof. }}
\newcommand{\qed}{\nopagebreak $\qquad$ $\square$ \vskip5pt}
\newcommand{\separate}{\vskip5pt}

\newcommand{\im}{\operatorname{Im}}

\newcommand{\tr}{\operatorname{Tr}}

\newcommand{\B}{\overline}
\newcommand{\HC}{\BB H_{\BB C}}
\newcommand{\HR}{\BB H_{\BB R}}

\newcommand{\M}{\operatorname{Mx}}

\begin{document}

\title{\bf An Invitation to Split Quaternionic Analysis}
\author{Matvei Libine}
\maketitle

\section{Introduction}

Six years after William Rowan Hamilton's discovery of quaternions,
in 1849 James Cockle introduced the algebra of split quaternions \cite{Co}.
(He called them ``coquaternions.'')
One way to define the split quaternions $\HR$ is by taking the standard
generators for the algebra of quaternions
$\BB H = \BB R 1 \oplus \BB R i \oplus \BB R j \oplus \BB R k$
and replacing the $i$ and $j$ with $\tilde\iota = \sqrt{-1} i$ and
$\tilde\jmath = -\sqrt{-1}j$ respectively, so that $\HR =
\BB R 1 \oplus \BB R \tilde\iota \oplus \BB R \tilde\jmath \oplus \BB R k$.
Another way to realize split quaternions is as real $2 \times 2$ matrices.
And yet another realization is
$$
\HR = \biggl\{
\begin{pmatrix} z_{11} & z_{12} \\ z_{21} & z_{22} \end{pmatrix} ; \:
z_{11}, z_{12}, z_{21}, z_{22} \in \BB C ,\:
z_{22} = \overline{z_{11}}, \: z_{21} = \overline{z_{12}} \biggr\}.
$$

I.~Frenkel initiated development of quaternionic analysis from
the point of view of representation theory of the conformal group
$SL(2,\BB H)$ and its Lie algebra $\mathfrak{sl}(2,\BB H)$. This approach
has already been proven very fruitful and in our joint work \cite{FL} we
push further the parallel with complex analysis and develop a rich theory.
In particular we show that the quaternionic analogue of
the Cauchy integral formula for the second order pole
$$
f'(w) = \frac 1{2\pi i} \oint \frac {f(z)\,dz}{(z-w)^2}
$$
are the differential operator 
$$
\M f = \nabla f \nabla - \square f^+
$$
defined on all holomorphic functions of four complex variables $f: \HC \to \HC$
and its integral presentation (note the square of the Fueter kernel)
$$
(\M f)(W) = \int_{C_4} \frac {(Z-W)^{-1}}{\det(Z-W)} \cdot
f(Z) \cdot \frac {(Z-W)^{-1}}{\det(Z-W)} \,dZ^4,
$$
where $dZ^4$ is the volume form, $C_4$ is a four cycle
homologous to $U(2) = \{Z \in \HC ;\: Z^*Z=1 \}$ sitting in the
complexified quaternionic space $\HC = \BB C \otimes \BB H$.
Since the constant functions on $\BB C$ are the holomorphic functions
annihilated by the operator $d/dz: f(z) \mapsto f'(z)$, their quaternionic
analogue is the kernel of the operator $\M$, which turns out to be the space
of solutions of a Euclidean version of the Maxwell equation for the gauge
potential.

We also identify the Feynman integrals associated to the diagrams

\vskip7pt

\centerline{\includegraphics{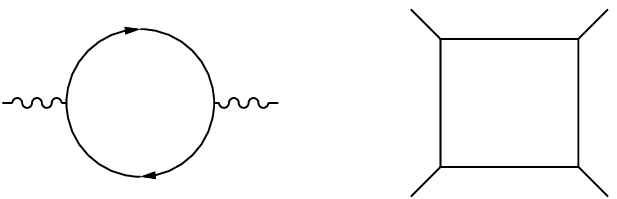}}
\nopagebreak
\centerline{Feynman diagrams}

\noindent
with the intertwining operators projecting certain natural unitary
representations of $\mathfrak{su}(2,2)$ onto their first irreducible
components. Then we conjecture that the other Feynman integrals also
admit an interpretation via quaternionic analysis and representation theory
as the projectors onto the other irreducible components.

With this representation theoretic approach it
quickly becomes evident that one has to consider the complexifications
$\HC = \BB C \otimes \BB H$ and $\mathfrak{sl}(4,\BB C)$ of $\BB H$ and
$\mathfrak{sl}(2,\BB H)$ and their real forms, such as the Minkowski space
$\BB M$ and $\mathfrak{su}(2,2)$ or the split quaternions $\HR$ and
$\mathfrak{sl}(4,\BB R)$. An important aspect of quaternionic analysis
is its ability to compare representation theories of various real forms,
and thus produce new results and make previously known results more explicit.

Just as (classical) quaternionic analysis is intimately
related to the representation theory of $SU(2)$, split
quaternionic analysis is related to the representation theory of
$SL(2,\BB R)$. The representation theory of
$SL(2,\BB R)$ is much richer than that of $SU(2)$
and exhibits most aspects of representations of higher rank real semisimple
Lie groups.
In particular, the group $SL(2,\BB R)$ exhibits a subtle aspect of
representation theory such as the separation of the
discrete and continuous series of unitary representations:
\begin{equation}  \label{separation}
L^2 \bigl( SL(2,\BB R) \bigr) \simeq
L^2_{discr} \bigl( SL(2,\BB R) \bigr) \oplus L^2_{cont} \bigl( SL(2,\BB R) \bigr).
\end{equation}
We study this decomposition from the quaternionic point of view.
The denominator of the Cauchy-Fueter kernel $\det(X-X_0)$ determines a
region in $\HC = \BB C \otimes \BB H$
$$
\Omega = \bigl\{ X_0 \in \HC ;\:
\det(X-X_0)\ne 0 \quad \forall X \in \HR \text{ with } \det(X)=1 \bigr\}.
$$
Loosely speaking, this region can be broken into several connected components,
two of which are open Ol'shanskii semigroups of $GL(2,\BB C)$.
Then the Cauchy-Fueter kernel $\frac{(X-X_0)^{-1}}{\det(X-X_0)}$ on $\HR$
can be expanded in terms of the $K$-types of the discrete or continuous
series of $SL(2,\BB R)$ depending on the choice of connected component of
$\Omega$ containing $X_0$.
It follows that the projectors onto the discrete and continuous series
of $SL(2,\BB R)$ can be expressed as
\begin{align*}
(\operatorname{P}^{discr}) f(X_0) &=
\frac 1 {2\pi^2} \int_{C_{discr}}
\frac{(X-X_0)^{-1}}{\det(X-X_0)} \cdot Dz \cdot f(X),  \\
(\operatorname{P}^{cont}) f(X_0) &=
\frac 1 {2\pi^2} \int_{C_{cont}}
\frac{(X-X_0)^{-1}}{\det(X-X_0)} \cdot Dz \cdot f(X),
\end{align*}
where $C_{discr}$ and $C_{cont}$ are certain three cycles in $\Omega$.
Note that these two integrals are identical to the Cauchy-Fueter formula,
except for the choice of contours of integration $C_{discr}$ and $C_{cont}$.
Moreover, by choosing appropriate cycles we can even get projectors onto
the holomorphic discrete series and antiholomorphic discrete series.
Such a geometric description of the decomposition (\ref{separation})
fits well into the Gelfand-Gindikin program initiated in \cite{GG}.
This relationship between quaternionic analysis and the separation of the
series for $SL(2, \BB R)$ will be the subject of our upcoming paper \cite{FL2}.


In view of how many connections to mathematical physics (classical)
quaternionic analysis has (see, for example, \cite{GT}), it is almost
guaranteed that split quaternionic analysis will have them as well.
At this point it is worth mentioning two very recent physics
papers \cite{ACCK} and \cite{MS} stressing the importance of the $(2,2)$
signature of the split quaternions $\HR$ as opposed to the traditional
Lorentzian signature of the Minkowski space $\BB M$.

Once split quaternionic analysis is sufficiently developed,
it can be generalized in many different ways.
Perhaps the most obvious direction is to extend the new results to
higher dimensions. There is a generalization of quaternionic analysis
known as Clifford analysis.
Let $\mathfrak{Cl}(V)$ be a Clifford algebra over a real finite-dimensional
vector space $V$ with a nondegenerate quadratic form $Q(x)$.
Then on $\BB R \oplus V$ we can introduce Dirac operators $D$ and
$D^+$ with coefficients in $\mathfrak{Cl}(V)$ so that $DD^+ = D^+D$ is
the wave operator on $\BB R \oplus V$ with symbol $x^2 -Q$.
We define a differentiable function $f: \BB R \oplus V \to \mathfrak{Cl}(V)$
to be (left) Clifford analytic (or monogenic) if $D^+f=0$.
(We think of it as an analogue of the Cauchy-Riemann equations.)
Slightly more generally, we can define Clifford analytic functions with
values in a $\mathfrak{Cl}(V)$-module.
When the quadratic form $Q$ is negative definite 
there is a well developed theory of Clifford analytic functions called
Clifford analysis (see, for example, \cite{BDS}, \cite{DSS} and \cite{GM}).
This theory generalizes (classical) quaternionic analysis.
Note that in this case $DD^+$ is the Laplacian on $\BB R \oplus V$, 
hence the components of Clifford analytic functions are harmonic.
Clifford analysis has deep connections with harmonic analysis, representations
of spin groups and index theory of Dirac operators.

Methods developed in split quaternionic analysis will extend to analysis
of Clifford analytic functions associated with quadratic forms of arbitrary
signature.
Thus we can consider split Clifford analysis as a ``real form'' of
complexified Clifford analysis, which in turn was introduced by Ryan in
\cite{R}.
On the other hand, the use of the wave equation for the study of harmonic
analysis on a hyperboloid goes back to Strichartz \cite{St}.
More recently, Kobayashi and {\O}rsted \cite{KO} study representations
of $O(p+1,q+1)$ in the space of solutions of the wave equation
$\square_{p,q}\phi=0$ on $\BB R^{p+q}$.
In this light it is natural to expect that this split Clifford analysis
will yield new results relating solutions of the wave equation,
representation theory and index theory of Dirac operators.
In particular, we expect to obtain concrete realizations of representations
of $O(p+1,q+1)$ in the space of solutions of Clifford analytic functions.
Since Clifford analysis is widely used by mathematical physicists,
it is very likely that they will find the split version at least
as useful as the classical one.


In this article we approach the split quaternions $\HR$ as a real form of
$\HC$, introduce the notion of regular functions and
give two different analogues of the Cauchy-Fueter formula valid for different
classes of functions.
This is done in parallel with (classical) quaternionic analysis.
We conclude the paper with an outline of our derivation of the projectors
$\operatorname{P}^{discr}$ and $\operatorname{P}^{cont}$ onto the discrete and
continuous series of $SL(2,\BB R)$.
Some contemporary reviews of quaternionic analysis are given in
\cite{Su} and \cite{CSSS}.

\section{The Quaternionic Spaces $\HC$, $\HR$ and $\BB M$}

In this article we use notations established in \cite{FL}.
In particular, $e_0$, $e_1$, $e_2$, $e_3$ denote the units of the classical
quaternions $\BB H$ corresponding to the more familiar $1$, $i$, $j$, $k$
(we reserve the symbol $i$ for $\sqrt{-1} \in \BB C$).
Thus $\BB H$ is an algebra over $\BB R$ generated by $e_0$, $e_1$, $e_2$, $e_3$,
and the multiplicative structure is determined by the rules
$$
e_0 e_i = e_i e_0 = e_i, \qquad
(e_i)^2 = - e_0, \qquad
e_ie_j=-e_ie_j, \qquad 1 \le i< j \le 3,
$$
and the fact that $\BB H$ is a division ring.
Next we consider the algebra of complexified quaternions
$\HC = \BB C \otimes \BB H$ (also known as biquaternions).
We define a complex conjugation on $\HC$ with respect to $\BB H$:
$$
Z = z^0e_0 + z^1e_1 + z^2e_2 + z^3e_3 \quad \mapsto \quad
Z^c = \B{z^0}e_0 + \B{z^1}e_1 + \B{z^2}e_2 + \B{z^3}e_3,
\qquad z^0,z^1,z^2,z^3 \in \BB C,
$$
so that $\BB H = \{ Z \in \HC ;\: Z^c=Z \}$.
The quaternionic conjugation on $\HC$ is defined by:
$$
Z = z^0e_0 + z^1e_1 + z^2e_2 + z^3e_3 \quad \mapsto \quad
Z^+ = z^0e_0 - z^1e_1 - z^2e_2 - z^3e_3,
\qquad z^0,z^1,z^2,z^3 \in \BB C;
$$
it is an anti-involution:
$$
(ZW)^+ = W^+ Z^+, \qquad \forall Z,W \in \HC.
$$
We will also use an involution
$$
Z \mapsto Z^- = - e_3 Z e_3 \qquad \text{(conjugation by $e_3$).}
$$
Then the complex conjugation, the quaternionic conjugation and the involution
$Z \mapsto Z^-$ commute with each other.

In this article we will be primarily interested in the space of
{\em split quaternions} $\HR$ which is a real form of $\HC$ defined by
$$
\HR = \{ Z \in \HC ;\: Z^{c-} = Z \}
= \{\text{$\BB R$ -span of
$e_0$, $\tilde e_1 = i e_1$, $\tilde e_2 = -i e_2$, $e_3$} \}.
$$
We will also consider the Minkowski space $\BB M$ which we regard as
another real form of $\HC$:
$$
\BB M = \{Z \in \HC ;\: Z^{c+} = - Z \},
$$
$\BB M$ is spanned over $\BB R$ by $\tilde e_0 = -i e_0$, $e_1$, $e_2$, $e_3$.

On $\HC$ we have a quadratic form $N$ defined by
$$
N(Z) = ZZ^+ = Z^+Z = (z^0)^2 + (z^1)^2 + (z^2)^2 + (z^3)^2,
$$
hence
$$
Z^{-1} = \frac {Z^+}{N(Z)}.
$$
The corresponding symmetric bilinear form on
$\HC$ is
\begin{equation}  \label{bilinear_form}
\langle Z, W \rangle = \frac 12 \tr(Z^+ W) = \frac 12 \tr(Z W^+),
\qquad Z,W \in \HC,
\end{equation}
where $\tr Z = 2z^0 = Z + Z^+$.
When this quadratic form is restricted to $\BB H$, $\HR$ and $\BB M$,
it has signature $(4,0)$, $(2,2)$ and $(3,1)$ respectively.
The real forms $\BB H$ and $\BB M$ have been studied in \cite{FL}, and
the signature $(1,3)$ is equivalent to $(3,1)$.
In this article we study the real form $\HR$ realizing the only
remaining signature $(2,2)$.

We will use the standard matrix realization of $\BB H$ so that
$$
e_0 = \begin{pmatrix} 1 & 0 \\ 0 & 1 \end{pmatrix}, \qquad
e_1 = \begin{pmatrix} 0 & -i \\ -i & 0 \end{pmatrix}, \qquad
e_2 = \begin{pmatrix} 0 & -1 \\ 1 & 0 \end{pmatrix}, \qquad
e_3 = \begin{pmatrix} -i & 0 \\ 0 & i \end{pmatrix},
$$
and
$$
\BB H = \{ Z \in \HC ;\: Z^c = Z \}
= \biggl\{
Z= \begin{pmatrix} z_{11} & z_{12} \\ z_{21} & z_{22} \end{pmatrix} \in \HC
; \: z_{22} = \overline{z_{11}}, \: z_{21} = - \overline{z_{12}} \biggr\}.
$$
Then $\HC$ can be identified with the algebra of all complex
$2 \times 2$ matrices:
$$
\HC = \biggl\{
Z= \begin{pmatrix} z_{11} & z_{12} \\ z_{21} & z_{22} \end{pmatrix}
; \: z_{ij} \in \BB C \biggr\},
$$
the quadratic form $N(Z)$ becomes $\det Z$ and the involution
$Z \mapsto Z^-$ becomes
$$
Z = \begin{pmatrix} z_{11} & z_{12} \\ z_{21} & z_{22} \end{pmatrix}
\mapsto Z^- =
\begin{pmatrix} 1 & 0 \\ 0 & -1 \end{pmatrix} Z
\begin{pmatrix} 1 & 0 \\ 0 & -1 \end{pmatrix}
= \begin{pmatrix} z_{11} & -z_{12} \\ -z_{21} & z_{22} \end{pmatrix}.
$$
The split quaternions $\HR$ and the Minkowski space $\BB M$ have matrix
realizations
$$
\HR = \biggl\{
Z= \begin{pmatrix} z_{11} & z_{12} \\ z_{21} & z_{22} \end{pmatrix} \in \HC
; \: z_{22} = \overline{z_{11}}, \: z_{21} = \overline{z_{12}} \biggr\}
$$
and
$$
\BB M = \biggl\{
Z= \begin{pmatrix} z_{11} & z_{12} \\ z_{21} & z_{22} \end{pmatrix} \in \HC
; \: z_{11}, z_{22} \in i \BB R, \: z_{21} = -\overline{z_{12}} \biggr\}.
$$
From this realization it is easy to see that the split quaternions form an
algebra over $\BB R$ isomorphic to $\mathfrak{gl} (2, \BB R)$,
the invertible elements in $\HR$, denoted by $\HR^{\times}$,
are nothing else but $GL(2,\BB R)$.
We regard $SL(2,\BB C)$ as a quadric $\{N(Z) = 1\}$ in $\HC$, and
we also regard $SU(1,1) \simeq SL(2,\BB R)$ as the set of real points of
this quadric:
\begin{multline}  \label{SU(1,1)}
SU(1,1) = \{ Z \in \HR ; \: N(Z) = 1 \}  \\
= \biggl\{ Z = \begin{pmatrix} z_{11} & z_{12} \\
\overline{z_{12}} & \overline{z_{11}} \end{pmatrix}
\in \HR ;\: \det Z = |z_{11}|^2 - |z_{12}|^2 = 1 \biggr\}.
\end{multline}

The algebra of split quaternions $\HR$ is spanned over $\BB R$
by the four matrices
$$
e_0 = \begin{pmatrix} 1 & 0 \\ 0 & 1 \end{pmatrix}, \qquad
\tilde e_1 = \begin{pmatrix} 0 & 1 \\ 1 & 0 \end{pmatrix}, \qquad
\tilde e_2 = \begin{pmatrix} 0 & i \\ -i & 0 \end{pmatrix}, \qquad
e_3 = \begin{pmatrix} -i & 0 \\ 0 & i \end{pmatrix},
$$
so
$$
\HR = \biggl\{ x^0 e_0 + x^1 \tilde e_1 + x^2 \tilde e_2 + x^3 e_3 =
\begin{pmatrix} x^0-ix^3 & x^1+ix^2 \\ x^1-ix^2 & x^0+ix^3 \end{pmatrix}
;\: x^0, x^1, x^2, x^3 \in \BB R \biggr\}.
$$
The quaternionic conjugation in this basis becomes
$$
e_0^+ = e_0, \quad \tilde e_1^+ = - \tilde e_1,
\quad \tilde e_2^+ = - \tilde e_2, \quad e_3^+ = -e_3.
$$
The multiplication rules for $\HR$ are:
\begin{center}
$e_0$ commutes with all elements of $\HR$, \\
$\tilde e_1$, $\tilde e_2$, $e_3$ anti-commute, \\
$e_0^2 = \tilde e_1^2 = \tilde e_2^2 = e_0$, \quad $e_3^2 = -e_0$,  \\
$\tilde e_1 \tilde e_2 = e_3$, \quad $\tilde e_2e_3 = -\tilde e_1$,
\quad $e_3 \tilde e_1= - \tilde e_2$.
\end{center}
The elements $e_0$, $\tilde e_1$, $\tilde e_2$, $e_3$
are orthogonal with respect to the bilinear form (\ref{bilinear_form})
and $\langle e_0, e_0 \rangle = \langle e_3, e_3 \rangle = 1$,
$\langle \tilde e_1, \tilde e_1 \rangle =
\langle \tilde e_2, \tilde e_2 \rangle = -1$.


The (classical) quaternions $\BB H$ are oriented so that
$\{e_0, e_1, e_2, e_3 \}$ is a positive basis.
Let $dV= dz^0 \wedge dz^1 \wedge dz^2 \wedge dz^3$ be the holomorphic 4-form
on $\HC$ determined by $dV(e_0,e_1,e_2,e_3)=1$, then the restriction
$dV \bigl |_{\BB H}$ is the Euclidean volume form corresponding to
$\{e_0, e_1, e_2, e_3 \}$.
On the other hand, the restriction $dV \bigl |_{\HR}$ is also real-valued
and hence determines an orientation on $\HR$ so that
$\{e_0, \tilde e_1, \tilde e_2, e_3 \}$ becomes a positively oriented basis.
Define a norm on $\HC$ by
$$
\|Z\| = \frac 1{\sqrt 2}
\sqrt{|z_{11}|^2 + |z_{12}|^2 + |z_{21}|^2 + |z_{22}|^2}, \qquad
Z = \begin{pmatrix} z_{11} & z_{12} \\ z_{21} & z_{22} \end{pmatrix} \in \HC,
$$
so that $\|e_i\|=1$, $0 \le i \le 3$.

In \cite{FL} we defined a holomorphic 3-form on $\HC$ with values in $\HC$
$$
Dz = e_0 dz^1 \wedge dz^2 \wedge dz^3 - e_1 dz^0 \wedge dz^2 \wedge dz^3
+ e_2 dz^0 \wedge dz^1 \wedge dz^3 - e_3 dz^0 \wedge dz^1 \wedge dz^2
$$
characterized by the property
\begin{multline*}
\langle Z_1, Dz(Z_2,Z_3,Z_4) \rangle =
\frac 12 \tr (Z_1^+, Dz (Z_2,Z_3,Z_4)) = dV(Z_1,Z_2,Z_3,Z_4), \\
\qquad \forall Z_1,Z_2,Z_3,Z_4 \in \HC.
\end{multline*}
Let $Dx = Dz \bigl |_{\HR}$ and $D \tilde x = Dz \bigl |_{\BB H}$.

\begin{prop}
The 3-form $Dx$ takes values in $\HR$.
If we write $X = x^0 e_0 + x^1 \tilde e_1 + x^2 \tilde e_2 + x^3 e_3 \in \HR$,
$x^0, x^1, x^2, x^3 \in \BB R$, then $Dx$ is given explicitly by
\begin{equation}  \label{Dx-explicit}
Dx = e_0 dx^1 \wedge dx^2 \wedge dx^3 + \tilde e_1 dx^0 \wedge dx^2 \wedge dx^3
- \tilde e_2 dx^0 \wedge dx^1 \wedge dx^3 - e_3 dx^0 \wedge dx^1 \wedge dx^2.
\end{equation}
\end{prop}

\begin{rem}
Clearly, the form $Dx$ satisfies the property
\begin{multline*}
\langle X_1, Dx(X_2,X_3,X_4) \rangle =
\frac 12 \tr (X_1^+, Dx (X_2,X_3,X_4)) = dV(X_1,X_2,X_3,X_4), \\
\qquad \forall X_1,X_2,X_3,X_4 \in \HR,
\end{multline*}
which could be used to define it.

It is also worth mentioning that in terms of the coordinates
$Z= \begin{pmatrix} z_{11} & z_{12} \\ z_{21} & z_{22} \end{pmatrix}$
on $\HC$, $z_{ij} \in \BB C$,
$$
Dz = \frac 12 \begin{pmatrix}
-dz_{11} \wedge dz_{12} \wedge dz_{21} &
-dz_{11} \wedge dz_{12} \wedge dz_{22}  \\
dz_{11} \wedge dz_{21} \wedge dz_{22} &
dz_{12} \wedge dz_{21} \wedge dz_{22}
\end{pmatrix},
$$
where we write $z_{ij} = x_{ij} + i y_{ij}$, $x_{ij}, y_{ij} \in \BB R$,
and $dz_{ij} = dx_{ij} + i dy_{ij}$.
\end{rem}

Let $U \subset \HR$ be an open region with piecewise smooth boundary
$\partial U$. We give a canonical orientation to $\partial U$ as follows.
The positive orientation of $U$ is determined by
$\{e_0, \tilde e_1, \tilde e_2, e_3 \}$.
Pick a smooth point $p \in \partial U$ and let $\overrightarrow{n_p}$
be a non-zero vector in $T_p\HR$ perpendicular to $T_p\partial U$ and
pointing outside of $U$.
Then $\{\overrightarrow{\tau_1}, \overrightarrow{\tau_2},
\overrightarrow{\tau_3}\} \subset T_p \partial U$ is positively oriented
in $\partial U$ if and only if
$\{\overrightarrow{n_p}, \overrightarrow{\tau_1}, \overrightarrow{\tau_2},
\overrightarrow{\tau_3}\}$ is positively oriented in $\HR$.

\begin{lem}  \label{restrictions}
Let $R \in \BB R$ be a constant, then we have the following restriction
formulas:
\begin{center}
\begin{tabular} {ccc}
$Dx \Bigl|_{\{X \in \HR ;\: N(X) = R\}} = \frac {X}{\|X\|} \,dS$, & \qquad
& $Dx \Bigl|_{\{X \in \HR ;\: \| X\| = R\}} = \frac {X^-}{\|X\|} \,dS
= \frac {X^-}R \,dS$,
\end{tabular}
\end{center}
where the sets $\{X \in \HR ;\: N(X) = R\}$ and
$\{X \in \HR ;\: \| X\| = R\}$ are oriented as boundaries of the open sets
 $\{X \in \HR ;\: N(X) < R\}$ and $\{X \in \HR ;\: \| X\| < R\}$
respectively, and $dS$ denotes the respective restrictions of the Euclidean
measure on $\HR$.
\end{lem}

\section{Regular Functions on $\BB H$ and $\HC$}

Recall that regular functions on $\BB H$ are defined using an analogue of the
Cauchy-Riemann equations. Write $\tilde X \in \BB H$ as
$\tilde X = \tilde x^0 e_0 + \tilde x^1 e_1 + \tilde x^2 e_2 + \tilde x^3 e_3$,
$\tilde x^0, \tilde x^1, \tilde x^2, \tilde x^3 \in \BB R$,
and factor the four-dimensional Laplacian operator $\square$ on $\BB H$
as a product of two Dirac operators
$$
\square = \frac {\partial^2}{(\partial \tilde x^0)^2} +
\frac {\partial^2}{(\partial \tilde x^1)^2} +
\frac {\partial^2}{(\partial \tilde x^2)^2} +
\frac {\partial^2}{(\partial \tilde x^3)^2}
= \nabla \nabla^+ = \nabla^+ \nabla,
$$
where
\begin{align*}
\nabla^+ &=
e_0 \frac{\partial}{\partial \tilde x^0}
+ e_1 \frac{\partial}{\partial \tilde x^1}
+ e_2 \frac{\partial}{\partial \tilde x^2}
+ e_3 \frac{\partial}{\partial \tilde x^3}
\qquad \text{and}  \\
\nabla &=
e_0 \frac{\partial}{\partial \tilde x^0}
- e_1 \frac{\partial}{\partial \tilde x^1}
- e_2 \frac{\partial}{\partial \tilde x^2}
- e_3 \frac{\partial}{\partial \tilde x^3}.
\end{align*}
The operators $\nabla^+$, $\nabla$ can be applied to functions
on the left and on the right.
For an open subset $U \subset \BB H$ and a
differentiable function $f$ on $U$ with values
in $\BB H$ or $\HC$, we say $f$ is {\em left-regular} if
$(\nabla^+ f)(\tilde X)=0$ for all $\tilde X \in U$,
and $f$ is {\em right-regular} if
$(f \nabla^+)(\tilde X)=0$ for all $\tilde X \in U$.

\begin{prop}
For any ${\cal C}^1$-function $f$ on $U \subset \BB H$
with values in $\BB H$ or $\HC$,
$$
d(f \cdot D\tilde x) =
df \wedge D\tilde x = (f \nabla^+) \,dV \bigr|_{\BB H},
\qquad
d(D\tilde x \cdot f) =
- D\tilde x \wedge df = (\nabla^+ f) \,dV \bigr|_{\BB H}.
$$
In particular,
$$
\nabla^+ f = 0 \quad \Longleftrightarrow \quad
D\tilde x \wedge df =0, \qquad
f \nabla^+ = 0 \quad \Longleftrightarrow \quad
df \wedge D\tilde x =0.
$$
\end{prop}

Following \cite{FL}, we say that a differential function
$f^{\BB C}: U^{\BB C} \to \HC$ defined on an open set $U^{\BB C} \subset \HC$ is
{\em holomorphic} if it is holomorphic with respect to the complex variables
$z^0, z^1, z^2, z^3$.
Then we define $f^{\BB C}$ to be {\em holomorphic left-regular} if
it is holomorphic and $\nabla^+ f^{\BB C} =0$.
Similarly, $f^{\BB C}$ is defined to be {\em holomorphic right-regular} if
it is holomorphic and $f^{\BB C} \nabla^+ =0$.

If we identify $\HC$ with complex $2 \times 2$ matrices
$\begin{pmatrix} z_{11} & z_{12} \\ z_{21} & z_{22} \end{pmatrix}$,
$z_{ij} \in \BB C$, then a function $f^{\BB C}: U^{\BB C} \to \HC$ is holomorphic
if and only if it is holomorphic with respect
to the complex variables $z_{ij}$, $1 \le i,j \le 2$.
Let us introduce holomorphic analogues of $\nabla^+$ and $\nabla$:
\begin{align*}
\nabla^+_{\BB C} &=
e_0 \frac{\partial}{\partial z^0} 
+ e_1 \frac{\partial}{\partial z^1}
+ e_2 \frac{\partial}{\partial z^2}
+ e_3 \frac{\partial}{\partial z^3}
\qquad \text{and}  \\
\nabla_{\BB C} &=
e_0 \frac{\partial}{\partial z^0} 
- e_1 \frac{\partial}{\partial z^1}
- e_2 \frac{\partial}{\partial z^2}
- e_3 \frac{\partial}{\partial z^3}.
\end{align*}
Then for a holomorphic function $f^{\BB C}: U^{\BB C} \to \HC$,
the following derivatives are equal:
$$
\nabla^+ f^{\BB C} = \nabla^+_{\BB C} f^{\BB C} =
2 \begin{pmatrix} \frac {\partial}{\partial z_{22}} &
- \frac {\partial}{\partial z_{21}}  \\
- \frac {\partial}{\partial z_{12}} &
\frac {\partial}{\partial z_{11}} \end{pmatrix} f^{\BB C},
\qquad
f^{\BB C} \nabla^+ = f^{\BB C} \nabla^+_{\BB C} =
2 f^{\BB C} \begin{pmatrix} \frac {\partial}{\partial z_{22}} &
- \frac {\partial}{\partial z_{21}}  \\
- \frac {\partial}{\partial z_{12}} &
\frac {\partial}{\partial z_{11}} \end{pmatrix},
$$
$$
\nabla f^{\BB C} = \nabla_{\BB C} f^{\BB C} =
2 \begin{pmatrix} \frac {\partial}{\partial z_{11}} &
\frac {\partial}{\partial z_{21}}  \\
\frac {\partial}{\partial z_{12}} &
\frac {\partial}{\partial z_{22}} \end{pmatrix} f^{\BB C},
\qquad
f^{\BB C} \nabla = f^{\BB C} \nabla_{\BB C} =
2 f^{\BB C} \begin{pmatrix} \frac {\partial}{\partial z_{11}} &
\frac {\partial}{\partial z_{21}}  \\
\frac {\partial}{\partial z_{12}} &
\frac {\partial}{\partial z_{22}} \end{pmatrix}.
$$

\begin{prop}
For any holomorphic function $f^{\BB C}: U^{\BB C} \to \HC$,
$$
\nabla^+_{\BB C} f^{\BB C} = 0 \quad \Longleftrightarrow \quad
Dz \wedge df^{\BB C} =0,
\qquad
f^{\BB C} \nabla^+_{\BB C}= 0 \quad \Longleftrightarrow \quad
df^{\BB C} \wedge Dz =0.
$$
\end{prop}

\begin{lem}  \label{closed}
We have:
\begin{enumerate}
\item
$\square \frac 1{N(Z)} = 0$;
\item
$\nabla_{\BB C} \frac 1{N(Z)}
= \frac 1{N(Z)} \nabla_{\BB C}
= -2 \frac {Z^{-1}}{N(Z)} = -2 \frac {Z^+} {N(Z)^2}$;
\item
$\frac {Z^{-1}}{N(Z)} = \frac {Z^+} {N(Z)^2}$
is a holomorphic left- and right-regular function defined wherever
$N(Z) \ne 0$;
\item
The form
$\frac {Z^{-1}}{N(Z)} \cdot Dz = \frac {Z^+} {N(Z)^2} \cdot Dz$ is a closed
holomorphic $\HC$-valued 3-form defined wherever $N(Z) \ne 0$.
\end{enumerate}
\end{lem}

\begin{lem}  \label{inverse-deriv}
Let $U^{\BB C} \subset \HC$ be an open subset.
For any differentiable function $F: U^{\BB C} \to \BB C$, we have:
$$
\nabla \bigl( F(Z^+) \bigr) = (\nabla^+ F)(Z^+),
\qquad
\nabla^+ \bigl( F(Z^+) \bigr) = (\nabla F)(Z^+),
$$
$$
\nabla \bigl( F(Z^{-1}) \bigr) =
- Z^{-1} \cdot (\nabla F)(Z^{-1}) \cdot Z^{-1}.
$$
\end{lem}

\section{Regular Functions on $\HR$}

We introduce linear differential operators on $\HR$
\begin{align*}
\nabla^+_{\BB R} &= e_0 \frac{\partial}{\partial x^0}
- \tilde e_1 \frac{\partial}{\partial x^1}
- \tilde e_2 \frac{\partial}{\partial x^2}
+ e_3 \frac{\partial}{\partial x^3}
\qquad \text{and}  \\
\nabla_{\BB R} &= e_0 \frac{\partial}{\partial x^0}
+ \tilde e_1 \frac{\partial}{\partial x^1}
+ \tilde e_2 \frac{\partial}{\partial x^2}
- e_3 \frac{\partial}{\partial x^3}
\end{align*}
which may be applied to functions on the left and on the right.

Fix an open subset $U \subset \HR$ and let $f$ be a
differentiable function on $U$ with values in $\HR$ or $\HC$.

\begin{df}
The function $f$ is {\em left-regular} if it satisfies
$$
(\nabla^+_{\BB R} f)(X) = e_0 \frac{\partial f}{\partial x^0}(X)
- \tilde e_1 \frac{\partial f}{\partial x^1}(X)
- \tilde e_2 \frac{\partial f}{\partial x^2}(X)
+ e_3 \frac{\partial f}{\partial x^3}(X) =0, \qquad \forall X \in U.
$$
Similarly, $f$ is {\em right-regular} if
$$
(f\nabla^+_{\BB R} )(X) = \frac{\partial f}{\partial x^0}(X)e_0
- \frac{\partial f}{\partial x^1}(X) \tilde e_1
- \frac{\partial f}{\partial x^2}(X) \tilde e_2
+ \frac{\partial f}{\partial x^3}(X)e_3 =0, \qquad \forall X \in U.
$$
\end{df}

We denote by $\square_{2,2}$ the ultrahyperbolic wave operator on $\HR$
which can be factored as follows:
$$
\square_{2,2} = \frac {\partial^2}{(\partial x^0)^2} -
\frac {\partial^2}{(\partial x^1)^2} -
\frac {\partial^2}{(\partial x^2)^2} +
\frac {\partial^2}{(\partial x^3)^2}
=\nabla_{\BB R} \nabla^+_{\BB R} = \nabla^+_{\BB R} \nabla_{\BB R}.
$$ 

\begin{prop}  \label{d(fDzg)}
For any ${\cal C}^1$-function $f: U \to \HR$ (or $f: U \to \HC$),
$$
d(f \cdot Dx) = df \wedge Dx = (f \nabla^+_{\BB R}) \,dV,
\qquad
d(Dx \cdot f) = - Dx \wedge df = (\nabla^+_{\BB R} f) \,dV.
$$
In particular, 
$$
\nabla^+_{\BB R} f = 0 \quad \Longleftrightarrow \quad Dx \wedge df =0,
\qquad
f \nabla^+_{\BB R} = 0 \quad \Longleftrightarrow \quad df \wedge Dx =0.
$$
\end{prop}

Let $U^{\BB C} \subset \HC$ be an open set. The restriction relations
$$
Dz \bigl |_{\HR} = Dx,
\qquad
Dz \bigl |_{\BB H} = D\tilde x
$$
imply that the restriction of a holomorphic left- or right-regular function
to $U^{\BB R} = U^{\BB C} \cap \HR$ produces a left- or right-regular function
on $U^{\BB R}$ respectively.
And the restriction of a holomorphic left- or right-regular function
to $U_{\BB H} = U^{\BB C} \cap \BB H$ also produces a left- or right-regular
function on $U_{\BB H}$ respectively.
Conversely, if one starts with, say, a left-regular function on $\HR$,
extends it holomorphically to a left-regular function on $\HC$
and then restricts the extension to $\BB H$, the resulting function
is left-regular on $\BB H$.

These properties of Dirac operators $\nabla^+_{\BB R}$ and $\nabla_{\BB R}$
and the notion of regular functions on $\HR$ are
in complete parallel with the Dirac operators $\nabla^+_{\BB M}$ and
$\nabla_{\BB M}$ and the corresponding notion of regular functions on $\BB M$
introduced in Section 3.2 in \cite{FL}.

\begin{prop}
Let $f^{\BB C}: U^{\BB C} \to \HC$ be a holomorphic function. Then
$$
\nabla^+ f^{\BB C} = \nabla^+_{\BB R} f^{\BB C} = \nabla^+_{\BB M} f^{\BB C}
= \nabla^+ f^{\BB C}, \qquad
f^{\BB C} \nabla^+ = f^{\BB C} \nabla^+_{\BB R} = f^{\BB C} \nabla^+_{\BB M}
= f^{\BB C} \nabla^+,
$$
$$
\nabla f^{\BB C} = \nabla_{\BB R} f^{\BB C} = \nabla_{\BB M} f^{\BB C}
= \nabla f^{\BB C}, \qquad
f^{\BB C} \nabla = f^{\BB C} \nabla_{\BB R} = f^{\BB C} \nabla_{\BB M} = f^{\BB C} \nabla.
$$
\end{prop}

Thus, essentially by design, the Dirac operators
$\nabla^+$, $\nabla^+_{\BB R}$, $\nabla^+_{\BB M}$, $\nabla^+_{\BB C}$ and
$\nabla$, $\nabla_{\BB R}$, $\nabla_{\BB M}$, $\nabla_{\BB C}$
(and hence the notions of regular functions on $\BB H$, $\HR$,
$\BB M$ and holomorphic regular functions on $\HC$) are all compatible.

\section{Fueter Formula for Holomorphic Regular Functions on $\HR$}

We are interested in extensions of the Cauchy-Fueter formula to functions
on $\HR$.
First we recall the classical version of the integral formula due to Fueter:

\begin{thm} [Cauchy-Fueter Formula \cite{F1, F2}]  \label{Fueter}
Let $U_{\BB H} \subset \BB H$ be an open bounded subset with piecewise
${\cal C}^1$ boundary $\partial U_{\BB H}$.
Suppose that $f(\tilde X)$ is left-regular on a
neighborhood of the closure $\overline{U_{\BB H}}$, then
$$
\frac 1 {2\pi^2} \int_{\partial U_{\BB H}}
\frac {(\tilde X - \tilde X_0)^{-1}}{N(\tilde X - \tilde X_0)}
\cdot D\tilde x \cdot f(\tilde X) =
\begin{cases}
f(\tilde X_0) & \text{if $\tilde X_0 \in U_{\BB H}$;}\\
0 & \text{if $\tilde X_0 \notin \overline{U_{\BB H}}$.}
\end{cases}
$$
If $g(\tilde X)$ is right-regular on a neighborhood of the
closure $\overline{U_{\BB H}}$, then
$$
\frac 1 {2\pi^2} \int_{\partial U_{\BB H}}
g(\tilde X) \cdot D \tilde x \cdot
\frac {(\tilde X - \tilde X_0)^{-1}}{N(\tilde X - \tilde X_0)} =
\begin{cases}
g(\tilde X_0) & \text{if $\tilde X_0 \in U_{\BB H}$;}\\
0 & \text{if $\tilde X_0 \notin \overline{U_{\BB H}}$.}
\end{cases}
$$
\end{thm}

Let $U \subset \HR$ be an open subset, and let $f$ be a ${\cal C}^1$-function
defined on a neighborhood of $\overline{U}$ such that $\nabla_{\BB R}^+ f =0$.
In this subsection we extend the Cauchy-Fueter integral formula to left-regular
functions which can be extended holomorphically to a
neighborhood of $\overline{U}$ in $\HC$.
(In other words, we assume that $f$ is a real-analytic function on
$\overline{U}$.)
Observe that the expression in the integral formula
$\frac {(\tilde X - \tilde X_0)^{-1}}{N(\tilde X - \tilde X_0)}
\cdot D\tilde x$ is nothing else but the restriction to $\BB H$ of the
holomorphic 3-form $\frac {(Z-\tilde X_0)^{-1}} {N(Z-\tilde X_0)} \cdot Dz$
which is the form from Lemma \ref{closed} translated by $\tilde X_0$.
For this reason we expect an integral formula of the kind
$$
f(X_0) = \frac 1 {2\pi^2} \int_{\partial U}
\biggl( \frac {(Z-X_0)^{-1}} {N(Z-X_0)} \cdot Dz \biggr)
\biggl|_{\HR} \cdot f(X),
\qquad \forall X_0 \in U.
$$
However, the integrand is singular wherever $N(Z-X_0)=0$.
We resolve this difficulty by deforming the contour of integration $\partial U$
in such a way that the integral is no longer singular.

Fix an $\epsilon \in \BB R$ and define an $\epsilon$-deformation
$h_{\epsilon}: \HC \to \HC$, $Z \mapsto Z_{\epsilon}$, by:
\begin{center}
\begin{tabular}{lccl}
$z_{11} \quad \mapsto \quad z_{11} + i\epsilon z_{11}$, & \qquad & \qquad &
$z_{12} \quad \mapsto \quad z_{12} - i\epsilon z_{12}$, \\
$z_{21} \quad \mapsto \quad z_{21} - i\epsilon z_{21}$, & \qquad & \qquad &
$z_{22} \quad \mapsto \quad z_{22} + i\epsilon z_{22}$.
\end{tabular}
\end{center}
Define a quadratic form on $\HC$
$$
S(Z) = z_{11}z_{22} + z_{12}z_{21}.
$$

\begin{lem}
We have the following identities:
$$
Z_{\epsilon} = Z + i \epsilon Z^-, \qquad
(Z_{\epsilon})^+ = Z^+ + i \epsilon Z^{+-},
$$
$$
N(Z_{\epsilon}) = (1-\epsilon^2) N(Z) + 2i\epsilon S(Z), \qquad
S(X) = \|X\|^2, \qquad \forall X \in \HR.
$$
\end{lem}

For $Z_0 \in \HC$ fixed, we use a notation
$$
h_{\epsilon, Z_0} (Z) = Z_0 + h_{\epsilon}(Z-Z_0) = Z + i\epsilon (Z-Z_0)^-.
$$

\begin{thm}  \label{holomorphic_Fueter}
Let $U \subset \HR$ be an open bounded subset with piecewise
${\cal C}^1$ boundary $\partial U$,
and let $f(X)$ be a ${\cal C}^1$-function defined on a
neighborhood of the closure $\overline{U}$ such that $\nabla^+_{\BB R} f =0$.
Suppose that $f$ extends to a holomorphic left-regular function
$f^{\BB C} : V^{\BB C} \to \HC$ with $V^{\BB C} \subset \HC$ an open 
subset containing $\overline{U}$, then
$$
- \frac 1 {2\pi^2} \int_{(h_{\epsilon, X_0})_*(\partial U)}
\frac {(Z-X_0)^{-1}} {N(Z-X_0)} \cdot Dz \cdot f^{\BB C}(Z) =
\begin{cases}
f(X_0) & \text{if $X_0 \in U$;}\\
0 & \text{if $X_0 \notin \overline{U}$,}
\end{cases}
$$
for all $\epsilon \ne 0$ sufficiently close to 0.
\end{thm}

\begin{rem}
For all $\epsilon \ne 0$ sufficiently close to 0 the contour of integration
$(h_{\epsilon, X_0})_*(\partial U)$ lies inside $V^{\BB C}$ and the integrand
is non-singular, thus the integrals are well-defined.
Moreover, we will see that the value of the integral becomes constant when
the parameter $\epsilon$ is sufficiently close to 0.
Of course, there is a similar formula for right-regular functions on $\HR$
as well.
\end{rem}

\pf
Let $M = \sup_{X \in \partial U} \|X-X_0\|$.
Without loss of generality we may assume that $V^{\BB C}$ is the
$\delta$-neighborhood of $\overline{U}$ for some $\delta >0$.
We will show that the integral formula holds for
$0 < |\epsilon| < \delta/M$.
Clearly, for this choice of $\epsilon$ the contour of integration
$(h_{\epsilon,X_0})_*(\partial U)$ lies inside $V^{\BB C}$ and, since the
integrand is a closed form, the integral stays constant
for $-\delta/M < \epsilon < 0$ and $0 < \epsilon < \delta/M$
(a priori the values of the integral may be different on these two intervals).

Let $S_r = \{ X \in \HR ;\: \|X-X_0\|^2 = r^2 \}$
and $B_r = \{ X \in \HR ;\: \|X-X_0\|^2 \le r^2 \}$
be the sphere and the closed ball of radius $r$ centered at $X_0$,
and choose $r>0$ sufficiently small so that $B_r \subset U$ and $r < \delta$.

\begin{lem}
Let $\tilde S_r = \{ \tilde X \in \BB H+X_0 ;\: \|\tilde X-X_0\|^2 = r^2 \}$
be the sphere oriented as the boundary of the open ball, then
$$
(h_{\epsilon, X_0})_*(\partial U) \quad \sim \quad
\begin{cases}
- \tilde S_r & \text{if $X_0 \in U$;} \\
0 & \text{if $X_0 \notin \overline{U}$}
\end{cases}
$$
as homology 3-cycles inside $\{Z \in V^{\BB C};\: N(Z - X_0) \ne 0 \}$.
\end{lem}

\pf
We give a proof for $\epsilon>0$; the case $\epsilon<0$ is similar.
As homology cycles in $\{Z \in V^{\BB C};\: N(Z - X_0) \ne 0 \}$,
$$
(h_{\epsilon, X_0})_*(\partial U) \quad \sim \quad
(h_{\epsilon, X_0})_*(S_r) \quad \sim \quad 
(h_{1, X_0})_*(S_r).
$$
If $X_0 \notin \overline{U}$, the cycle $(h_{1, X_0})_*(S_r)$ is homologous to
zero.

Assume now $X_0 \in U$.
Let $\P_{\BB H}$ be the projection $\HC \twoheadrightarrow \BB H$ defined
by 
\begin{multline*}
Z = (\tilde x^0 + i \tilde y^0) e_0
+ (\tilde x^1 + i \tilde y^1) e_1
+ (\tilde x^2 + i \tilde y^2) e_2
+ (\tilde x^3 + i \tilde y^3) e_3  \\
\mapsto
\tilde X = \tilde x^0 e_0 + \tilde x^1 e_1
+ \tilde x^2 e_2 + \tilde x^3 e_3,
\qquad
\tilde x^0,  \tilde x^1,  \tilde x^2, \tilde x^3,
\tilde y^0,  \tilde y^1,  \tilde y^2,  \tilde y^3 \in \BB R,
\end{multline*}
and let $\P_{\BB H+X_0} : \HC \twoheadrightarrow \BB H + X_0$,
be the projection
$\P_{\BB H+X_0} (Z) = \P_{\BB H} (Z-X_0) + X_0$.
We describe the supports of the cycles involved:
\begin{align*}
|S_r| &=
\{ X_0 + ae_0 + b \tilde e_1 + c \tilde e_2 + de_3
;\: a^2 + b^2 + c^2 + d^2 = r^2 \},  \\
|(h_{1, X_0})_*(S_r)| &= \{ X_0 +
(1+i)ae_0 + (1-i)b \tilde e_1 + (1-i)c \tilde e_2 + (1+i)de_3;  \\
&\hspace{3in} a^2 + b^2 + c^2 + d^2 = r^2 \},  \\
|(\P_{\BB H + X_0} \circ h_{1,X_0})_*(S_r)| &=
\{ X_0 + ae_0 + be_1 - ce_2 + de_3 ;\:  a^2 + b^2 + c^2 + d^2 = r^2 \}
= |\tilde S_r|.
\end{align*}
Moreover,
$$
(\P_{\BB H + X_0} \circ h_{1,X_0})_*(S_r) = - \tilde S_r
$$
as homology cycles.
It is easy to see that this projection provides a homotopy between
$(h_{1, X_0})_*(S_r)$ and $-S_r$, hence the lemma.
\qed

By Stokes'
$$
\int_{(h_{\epsilon, X_0})_*(\partial U)}
\frac {(Z-X_0)^{-1}} {N(Z-X_0)} \cdot Dz \cdot f^{\BB C}(Z)
=
\int_{-\tilde S_r}
\frac {(Z-X_0)^{-1}} {N(Z-X_0)} \cdot Dz \cdot f^{\BB C}(Z)
\qquad \text{if $X_0 \in U$,}
$$
and zero if $X_0 \notin \overline{U}$.
Finally, by the Fueter formula for the usual quaternions
(Theorem \ref{Fueter}), the last integral is $-2\pi^2 f(X_0)$.
(Alternatively, one can let $r \to 0^+$ and show directly that the integral
remains unchanged and at the same time approaches $-2\pi^2 f(X_0)$ in the same
way the Cauchy and Cauchy-Fueter formulas are proved.)
\qed

For a Cauchy-Fueter formula for regular functions on $\BB M$
that extend to holomorphic regular functions on $\HC$
see Section 3.3 in \cite{FL}.

\section{Fueter Formula for Regular Functions on $\HR$}


In this section we prove an analogue of the Cauchy-Fueter formula
for smooth left-regular functions on $\HR$ which are not necessarily real
analytic and do not necessarily have holomorphic extensions.
As a ``trade-off'' for working with ``bad'' functions the proofs become
much more involved.

\begin{thm}  \label{split_Fueter}
Let $U \subset \HR$ be a bounded open region with smooth boundary
$\partial U$.
Let $f: U \to \HC$ be a function which extends to a real-differentiable
function on an open neighborhood $V \subset \HR$ of the closure
$\overline{U}$ such that $\nabla^+_{\BB R} f = 0$.
Then, for any point $X_0 \in \HR$ such that $\partial U$ intersects the cone
$\{ X \in \HR ;\: N(X-X_0) =0 \}$ transversally, we have:
$$
\lim_{\epsilon \to 0} \frac {-1}{2\pi^2} \int_{\partial U}
\frac {(X-X_0)^+} {\bigl( N(X-X_0) +i\epsilon \|X-X_0\|^2 \bigr)^2}
\cdot Dz \cdot f(X) =
\begin{cases}
f(X_0) & \text{if $X_0 \in U$;}\\
0 & \text{if $X_0 \notin \overline{U}$.}
\end{cases}
$$
\end{thm}

\pf
The case $X_0 \notin \overline{U}$ is easier, so we assume $X_0 \in U$.
Using Proposition \ref{d(fDzg)}, we get
\begin{multline}  \label{d()}
d \biggl( \frac {(X-X_0)^+} {\bigl( N(X-X_0) +i\epsilon \|X-X_0\|^2 \bigr)^2}
\cdot Dz \cdot f(X) \biggr)  \\
=
\biggl( \frac {(X-X_0)^+} {\bigl( N(X-X_0) +i\epsilon \|X-X_0\|^2 \bigr)^2}
\nabla_{\BB R}^+ \biggr) f(X) \,dV  \\
= 4i\epsilon \frac {\|X-X_0\|^2 - (X-X_0)^+(X-X_0)^-}
{\bigl( N(X-X_0) +i\epsilon \|X-X_0\|^2 \bigr)^3} f(X) \,dV.
\end{multline}
In particular, expression (\ref{d()}) tends to zero pointwise when
$\epsilon \to 0$ except for those $X$ which lie on the cone
$\{ X \in \HR ;\: N(X-X_0) =0 \}$, and we need to be very careful there.
By translation we can assume that $X_0=0$.
Let $S_{r} = \{ X \in \HR ;\: \|X\|^2 = r^2 \}$
and $B_{r} = \{ X \in \HR ;\: \|X\|^2 \le r^2 \}$
be the sphere and the closed ball of radius $r$.
By Stokes'
\begin{multline*}
\int_{\partial U} \frac {X^+} {\bigl( N(X) +i\epsilon \|X\|^2 \bigr)^2}
\cdot Dz \cdot f(X)
= \int_{U \setminus  B_{r}}
d \biggl( \frac {X^+} {\bigl( N(X) +i\epsilon \|X\|^2 \bigr)^2}
\cdot Dz \cdot f(X) \biggr)  \\
+ \int_{S_{r}} \frac {X^+} {\bigl( N(X) +i\epsilon r^2 \bigr)^2}
\cdot Dz \cdot f(X).
\end{multline*}
We will show that, as $\epsilon \to 0$, the first integral on the right hand
side tends to zero (this part is non-trivial and uses that the cone
$\{ N(X-X_0) =0 \}$ intersects $\partial U$ transversally).
On the other hand, as $\epsilon \to 0$ and $r \to 0^+$, the second integral
tends to $-2\pi^2 f(0)$.
The proof is essentially a series of integration by parts computations.
In regular coordinates $(x^0,x^1,x^2,x^3)$ we have
$X = x^0e_0 + x^1 \tilde e_1 + x^2 \tilde e_2 + x^3e_3$,
and $N(X) = (x^0)^2 - (x^1)^2 - (x^2)^2 + (x^3)^2$.
For computing purposes we replace $(x^0,x^1,x^2,x^3)$ with spherical coordinates
$(\rho,\theta, \phi, \psi)$ so that
\begin{equation}  \label{spherical_coords}
\begin{matrix}
x^0 = \rho \cos \theta \cos \phi  \\
x^1 = \rho \sin \theta \sin \psi  \\
x^2 = \rho \sin \theta \cos \psi  \\
x^3 = \rho \cos \theta \sin \phi
\end{matrix}
\qquad
\begin{matrix}
\rho \ge 0,  \\
0 \le \theta \le \pi/2,  \\
0 \le \phi \le 2\pi,  \\
0 \le \psi \le 2\pi.
\end{matrix}
\end{equation}
This is an orientation-preserving change of coordinates, and the vector fields
$\{ \frac{\partial}{\partial \theta}, \frac{\partial}{\partial \phi},
\frac{\partial}{\partial \psi} \}$ form a positively-oriented frame on $S_r$.
Then
$$
N(X) = \rho^2 \cos(2\theta),
\qquad N(X) + i\epsilon \|X\|^2
= \rho^2 \bigl( \cos(2\theta) + i\epsilon \bigr),
$$
and the equation $N(X) = 0$ becomes $\theta = \pi/4$.

Recall that the function $\frac 1{x^n}$ which is singular at $x=0$ can
be regularized as a distribution in two different ways,
$\frac 1{(x + i0)^n}$ and $\frac 1{(x - i0)^n}$,
so that a test function $g(x)$ is being sent into
$\bigl\langle \frac 1{(x+i0)^n}, g(x) \bigr\rangle =
\lim_{\epsilon \to 0^+} \int \frac {g(x) \,dx}{(x+i\epsilon)^n}$
or $\bigl\langle \frac 1{(x-i0)^n}, g(x) \bigr\rangle =
\lim_{\epsilon \to 0^-} \int \frac {g(x) \,dx}{(x+i\epsilon)^n}$.
By a similar fashion we have the following lemma:

\begin{lem}  \label{regularization}
Fix a $\theta_0 \in (0, \frac {\pi}4)$,
and let $n$ be a positive integer, then we have two distributions
which send a test function $g(\theta)$ into the limits
$$
\lim_{\epsilon \to 0^+}
\int_{\frac{\pi}4-\theta_0}^{\frac{\pi}4+\theta_0} \frac {g(\theta) \, d\theta}
{\bigl( \cos (2\theta) + i\epsilon \bigr)^n}
\qquad \text{and} \qquad
\lim_{\epsilon \to 0^-}
\int_{\frac{\pi}4-\theta_0}^{\frac{\pi}4+\theta_0} \frac {g(\theta) \, d\theta}
{\bigl( \cos (2\theta) + i\epsilon \bigr)^n}.
$$
\end{lem}

\pf
We need to show that the limits exist and depend continuously
on the test function $g(\theta)$.
We do it by induction on $n$ using integration by parts.
If $n=1$,
\begin{multline*}
\int_{\frac{\pi}4-\theta_0}^{\frac{\pi}4+\theta_0}
\frac {g(\theta) \, d\theta}{\cos (2\theta) +i\epsilon}
=
\int_{\frac{\pi}4-\theta_0}^{\frac{\pi}4+\theta_0}
\frac {2\sin(2\theta)}{\cos (2\theta) +i\epsilon}
\cdot \frac {g(\theta)}{2\sin(2\theta)} \, d\theta  \\
=
\int_{\frac{\pi}4-\theta_0}^{\frac{\pi}4+\theta_0}
\log \bigl(\cos (2\theta) + i \epsilon \bigr) \cdot
\frac d{d\theta} \biggl( \frac {g(\theta)}{2\sin(2\theta)} \biggr) \,d\theta
-
\log \bigl( \cos (2\theta) + i \epsilon \bigr) \cdot
\frac {g(\theta)}{2\sin(2\theta)}
\biggl|_{\frac{\pi}4-\theta_0}^{\frac{\pi}4+\theta_0}.
\end{multline*}
For the purpose of this integration, the complex logarithm function is defined
on the complex plane $\BB C$ minus the negative real axis, and the values of
the logarithm lie in the strip $\{ z \in \BB C ;\: -\pi < \im z < \pi \}$.
The function $\log \bigl( \cos (2\theta) + i \epsilon \bigr)$ is
integrable for all values of $\epsilon$, including $\epsilon=0$, hence the
limits as $\epsilon \to 0^{\pm}$ exist and depend continuously on $g(\theta)$.

Now suppose that $n>1$, then
\begin{multline*}
\int_{\frac{\pi}4-\theta_0}^{\frac{\pi}4+\theta_0}
\frac {g(\theta) \, d\theta}{\bigl( \cos (2\theta) +i\epsilon \bigr)^n}
=
\int_{\frac{\pi}4-\theta_0}^{\frac{\pi}4+\theta_0}
\frac {2\sin(2\theta)} {\bigl( \cos (2\theta) +i\epsilon \bigr)^n}
\cdot \frac {g(\theta)}{2\sin(2\theta)} \, d\theta  \\
=
\frac 1{n-1} \frac 1{\bigl( \cos (2\theta) +i\epsilon \bigr)^{n-1}}
\cdot \frac {g(\theta)}{2\sin(2\theta)}
\biggl|_{\frac{\pi}4-\theta_0}^{\frac{\pi}4+\theta_0}  \\
- \frac 1{n-1}
\int_{\frac{\pi}4-\theta_0}^{\frac{\pi}4+\theta_0}
\frac 1{\bigl( \cos (2\theta) +i\epsilon \bigr)^{n-1}} \cdot
\frac d{d\theta}
\biggl( \frac {g(\theta)}{2\sin(2\theta)} \biggr) \, d\theta,
\end{multline*}
and the result follows by induction in $n$.
\qed

\begin{lem}
$$
\lim_{\epsilon \to 0} \int_{U \setminus B_{r}}
d \biggl( \frac {X^+} {\bigl( N(X) +i\epsilon \|X\|^2 \bigr)^2}
\cdot Dz \cdot f(X) \biggr) =0.
$$
\end{lem}

\pf
We have seen that
\begin{multline*}
\int_{U \setminus B_{r}} d \biggl(
\frac {X^+}{\bigl( N(X) +i\epsilon \|X\|^2 \bigr)^2}
\cdot Dz \cdot f(X) \biggr)  \\
=
4i\epsilon \int_{U \setminus B_{r}} \frac {\|X-X_0\|^2 - (X-X_0)^+(X-X_0)^-}
{\bigl( N(X) +i\epsilon \|X\|^2 \bigr)^3} f(X) \,dV.
\end{multline*}
Writing the right hand side integral in the spherical coordinates
(\ref{spherical_coords}) and integrating out the variables
$\rho$, $\phi$, $\psi$ we obtain an integral of the type
\begin{equation}  \label{integrals_to_zero}
4i\epsilon \int_0^{\frac {\pi}2} \frac {g(\theta) \,d\theta}
{\bigl( \cos (2\theta) +i\epsilon \bigr)^3},
\end{equation}
for some function $g(\theta)$.
By assumption, the boundary $\partial U$ is smooth,
compact and intersects the cone $\{N(X) =0\} = \{\theta=\pi/4\}$
transversally, hence the function $g(\theta)$ is smooth
at least for $\theta$ lying in some interval
$[\frac{\pi}4-\theta_0, \frac{\pi}4+\theta_0]$ with
$\theta_0 \in (0, \frac {\pi}4)$.
It follows from Lemma \ref{regularization} that the limit of
(\ref{integrals_to_zero}) as $\epsilon \to 0$ is zero.
\qed


\begin{lem}
$$
\int_{S_{r}} \frac {X^+} {\bigl( N(X) +i\epsilon r^2 \bigr)^2} \cdot Dz
=
r \int_{S_{r}} \frac {dS}{\bigl( N(X) + i\epsilon r^2 \bigr)^2}
=
-\frac {2\pi^2}{1 + \epsilon^2}.
$$
\end{lem}

\pf
From Lemma \ref{restrictions} we see that
$$
\frac {X^+} {\bigl( N(X) +i\epsilon r^2 \bigr)^2} \cdot Dz \biggl|_{S_{r}} =
\frac {X^+X^-} {\bigl( N(X) + i \epsilon r^2 \bigr)^2} \cdot \frac {dS}{r}.
$$
Notice that the involution $X \mapsto X^-$ preserves the sphere,
its orientation, its volume form $dS$ and replaces
$\frac {X^+X^-}{(N(X) + i\epsilon r^2)^2}$ with
$\frac {X^{+-} X}{(N(X) + i \epsilon r^2)^2}$.
Therefore, using $X^+ X^- + X^{+-} X = 2\|X\|^2$,
spherical coordinates (\ref{spherical_coords})
and $dS = r^3 \sin \theta \cos \theta\,d\psi d\phi d\theta$,
\begin{multline*}
\int_{S_{r}} \frac {X^+ X^-}{\bigl( N(X) + i\epsilon r^2 \bigr)^2}
\cdot \frac {dS}{r}
= \frac 12 \int_{S_{r}} \frac {X^+ X^- + X^{+-} X}
{\bigl( N(X) + i \epsilon r^2 \bigr)^2} \cdot \frac {dS}{r}
= \int_{S_{r}} \frac {r \,dS} {\bigl( N(X) + i\epsilon r^2 \bigr)^2}  \\
= \int_{\theta=0}^{\theta=\pi/2} \int_{\phi=0}^{\phi=2\pi}
\int_{\psi=0}^{\psi=2\pi}
\frac {\sin\theta \cos \theta \,d\psi d\phi d\theta}
{\bigl( \cos (2\theta) + i\epsilon \bigr)^2}
=
2\pi^2 \int_{\theta=0}^{\theta=\pi/2} \frac {\sin (2\theta) \,d\theta}
{\bigl( \cos (2\theta) + i\epsilon \bigr)^2}  \\
= \frac {\pi^2}{\cos (2\theta) + i\epsilon}
\biggl|_{\theta=0}^{\theta=\pi/2}
= \frac {\pi^2}{-1 + i\epsilon} - \frac {\pi^2}{1 + i\epsilon}
= - \frac{2\pi^2}{1 + \epsilon^2}.
\end{multline*} 
\qed

\begin{lem}
$$
\lim_{r \to 0^+} \biggl( \lim_{\epsilon \to 0}
\int_{S_{r}} \frac {X^+ \cdot Dz \cdot f(X)}{\bigl( N(X) +i\epsilon r^2 \bigr)^2}
\biggr) = -2\pi^2 f(0).
$$
\end{lem}

\pf
By Lemma \ref{restrictions} we have:
\begin{equation}  \label{epsilon_to_0}
\lim_{\epsilon \to 0}
\int_{S_{r}} \frac {X^+ \cdot Dz \cdot f(X)} {\bigl( N(X) +i\epsilon r^2 \bigr)^2}
= \lim_{\epsilon \to 0} \int_{S_{r}}
\frac {X^+ X^- \cdot f(X)} {\bigl( N(X) + i\epsilon r^2 \bigr)^2} \,\frac{dS}{r}
\end{equation}
If the function $f(X)$ were constant we would be done by previous lemma.
However, we cannot argue that since $r \to 0^+$, $f(X)$ is close to $f(0)$
and so may be treated like a constant because there can be derivatives of
$f(X)$ involved. So an integration by parts argument will be needed.

First we compute (writing $X = x^0e_0 + x^1 \tilde e_1 + x^2 \tilde e_2 + x^3e_3$)
$$
X^+ X^- =
\begin{pmatrix}
(x^0)^2+(x^1)^2+(x^2)^2+(x^3)^2 & 2(x^2x^3 - x^0x^1) - 2i (x^0x^2 + x^1x^3) \\
2(x^2x^3 - x^0x^1) + 2i (x^0x^2 + x^1x^3) & (x^0)^2+(x^1)^2+(x^2)^2+(x^3)^2
\end{pmatrix},
$$
and using the spherical coordinates (\ref{spherical_coords}) we can rewrite
$$
X^+ X^- =
\rho^2 \begin{pmatrix} 1 & 0 \\ 0 & 1 \end{pmatrix}
+ \rho^2 \sin(2\theta)
\sin(\psi-\phi) \begin{pmatrix} 0 & 1 \\ 1 & 0 \end{pmatrix}
+ \rho^2 \sin(2\theta)
\cos(\psi-\phi) \begin{pmatrix} 0 & -i \\ i & 0 \end{pmatrix}.
$$
Thus (\ref{epsilon_to_0}) can be rewritten as
$$
\lim_{\epsilon \to 0} \Biggl( r \int_{S_{r}}
\frac {f(X) \,dS} {\bigl( N(X) + i\epsilon r^2 \bigr)^2}
+
r \int_{S_{r}}
\frac { \sin(\psi-\phi) \begin{pmatrix} 0 & 1 \\ 1 & 0 \end{pmatrix}
+ \cos(\psi-\phi) \begin{pmatrix} 0 & -i \\ i & 0 \end{pmatrix}}
{\bigl( N(X) + i\epsilon r^2 \bigr)^2} \cdot \sin(2\theta) f \,dS \Biggr).
$$
We have:
\begin{multline*}
r \int_{S_{r}} \frac {f(X) \,dS} {\bigl( N(X) + i\epsilon r^2 \bigr)^2}
=
\frac 12 \int_{\theta=0}^{\theta=\pi/2} \int_{\phi=0}^{\phi=2\pi}
\int_{\psi=0}^{\psi=2\pi}
\frac {\sin (2\theta) f \,d\psi d\phi d\theta}
{\bigl( \cos (2\theta) + i\epsilon \bigr)^2}  \\
= \frac 14 \biggl(
\int_{\phi=0}^{\phi=2\pi} \int_{\psi=0}^{\psi=2\pi}
\frac {f}{\cos (2\theta) + i\epsilon}
\biggl|_{\theta=0}^{\theta=\pi/2} \,d\psi d\phi
- \int_{\theta=0}^{\theta=\pi/2} \int_{\phi=0}^{\phi=2\pi}
\int_{\psi=0}^{\psi=2\pi}
\frac { \frac {\partial f}{\partial \theta} \,d\psi d\phi d\theta}
{\cos (2\theta) + i\epsilon} \biggr).
\end{multline*}
By the chain rule $\frac {\partial f}{\partial \theta}$ can be written as
$\rho \cdot g(X)$ for some smooth function $g(X)$, thus the second
integral is $r$ times an expression from Lemma \ref{regularization}.
Taking limits $\epsilon \to 0$ and then $r \to 0^+$,
the first integral tends to $-2\pi^2 f(0)$ and the second tends to zero.
The second term inside the limit is
\begin{multline}  \label{integral_2}
r \int_{S_{r}}
\frac { \sin(\psi-\phi) \begin{pmatrix} 0 & 1 \\ 1 & 0 \end{pmatrix}
+ \cos(\psi-\phi) \begin{pmatrix} 0 & -i \\ i & 0 \end{pmatrix}}
{\bigl( N(X) + i\epsilon r^2 \bigr)^2} \cdot \sin(2\theta) f \,dS  \\
= \frac 12
\int_{\theta=0}^{\theta=\pi/2} \int_{\phi=0}^{\phi=2\pi} \int_{\psi=0}^{\psi=2\pi}
\frac {\sin(\psi-\phi) \begin{pmatrix} 0 & 1 \\ 1 & 0 \end{pmatrix}
+ \cos(\psi-\phi) \begin{pmatrix} 0 & -i \\ i & 0 \end{pmatrix}}
{\bigl( \cos (2\theta) + i\epsilon \bigr)^2}
\cdot \sin^2 (2\theta) f \,d\psi d\phi d\theta.
\end{multline}
But
\begin{align*}
\int_{\psi=0}^{\psi=2\pi} \sin (\psi-\phi)f \,d\psi
&= \int_{\psi=0}^{\psi=2\pi}
\cos (\psi-\phi) \frac {\partial f}{\partial \psi}f \,d\psi,  \\
\int_{\psi=0}^{\psi=2\pi} \cos (\psi-\phi)f \,d\psi
&= - \int_{\psi=0}^{\psi=2\pi}
\sin (\psi-\phi) \frac {\partial f}{\partial \psi}f \,d\psi.
\end{align*}
By the chain rule $\frac {\partial f}{\partial \psi}$ can be written as
$\rho \cdot h(X)$ for some smooth function $h(X)$, thus the right hand side
of (\ref{integral_2})
is $r$ times an expression from Lemma \ref{regularization}.
When we take limits first as $\epsilon \to 0$ and then as $r \to 0^+$,
integral (\ref{integral_2}) tends to zero.
\qed

This concludes our proof of Theorem \ref{split_Fueter}.
\qed

\section{Separation of the Series for $SL(2,\BB R)$}

What makes the representation theory of $SL(2,\BB R)$ more interesting
than that of $SU(2)$ is having the separation of the series into
discrete and continuous components.
Instead of $SL(2,\BB R)$ we prefer to work with $SU(1,1)$ sitting inside
$\HR$, as in (\ref{SU(1,1)}).
In this section we outline a relationship between split quaternionic analysis
and the decomposition
$$
L^2 \bigl( SU(1,1) \bigr) \simeq
L^2_{discr} \bigl( SU(1,1) \bigr) \oplus L^2_{cont} \bigl( SU(1,1) \bigr).
$$

The denominator of the Cauchy-Fueter kernel $N(X-X_0)$ determines a region
in $\HC$
$$
\Omega = \bigl\{ X_0 \in \HC ;\:
N(X-X_0)\ne 0 \quad \forall X \in \HR \text{ with } N(X)=1 \bigr\}.
$$
This region contains two open Ol'shanskii semigroups of $GL(2,\BB C)$
$$
\Gamma^0 = \{ Z \in \HC ;\: Z^* \tilde e_3 Z - \tilde e_3
\text{ is positive definite}\}, \qquad
$$
$$
\overline{\Gamma^0} = (\Gamma^0)^{-1} =
\{ Z \in \HC ;\: Z^* \tilde e_3 Z - \tilde e_3 \text{ is negative definite}\},
$$
where
$\tilde e_3 = ie_3 = \begin{pmatrix} 1 & 0 \\ 0 & -1 \end{pmatrix} \in \HC$.
Following \cite{KouO} we can show that over $\Gamma^0$ and $\overline{\Gamma^0}$
the Cauchy-Fueter kernel $\frac{(Z-X_0)^{-1}}{N(Z-X_0)}$ can be expanded in
terms of the $K$-types of the discrete series of $SU(1,1)$.
Thus we obtain the following integral formula for the projector onto the
discrete series of $SU(1,1))$:
\begin{equation} \label{discrete-proj}
(\operatorname{P}^{discr}) f(X_0) = \frac 1 {2\pi^2} \int_{C_{discr}}
\frac{(Z-X_0)^{-1}}{N(Z-X_0)} \cdot Dz \cdot f(Z),
\end{equation}
where $C_{discr}$ is a certain three cycle in $\Gamma^0 \cup \overline{\Gamma^0}$.
(Strictly speaking, this operator is not a projector because it has
eigenvalues $\pm 1$ on the discrete series subspace,
but its square does give a projection onto the discrete series.)
Furthermore, we can decompose the cycle $C_{discr}$ into a sum of two cycles
$C_{discr}^+ + C_{discr}^-$ with $C_{discr}^+$ and $C_{discr}^-$ lying in
$\overline{\Gamma^0}$ and $\Gamma^0$ respectively. Then integration over
$C_{discr}^+$ (respectively $C_{discr}^-$) produces a ``projector'' onto the
holomorphic (respectively antiholomorphic) discrete series.

Our next goal is to obtain a projector onto the continuous series component
$\operatorname{P}^{cont}$.
We expect that $\operatorname{P}^{cont}$ will be given by the same
formula (\ref{discrete-proj}) but with a different choice of the cycle
of integration, quite possibly supported in
$\Omega \setminus (\Gamma^0 \cup \overline{\Gamma^0})$.
To get $\operatorname{P}^{cont}$ we use a conformal map
$\gamma: \HC \to \HC$ which sends
$$
\HR \,\tilde\to\, \BB M \qquad \text{and} \qquad
SU(1,1) \,\tilde\to\, \widetilde{H}^3 \qquad \text{(with singularities),}
$$
where $\widetilde{H}^3$ is the unit hyperboloid of one sheet in $\BB M$;
we call $\gamma$ the ``Cayley transform''.
The hyperboloid $\widetilde{H}^3$ can be identified with $SL(2,\BB C)/SU(1,1)$
and is usually called the imaginary Lobachevski space.
The group $SL(2,\BB C)$ acts naturally on $L^2(\widetilde{H}^3)$ and
decomposes into the discrete and continuous components:
$$
L^2(\widetilde{H}^3) \simeq
L^2_{discr}(\widetilde{H}^3) \oplus L^2_{cont}(\widetilde{H}^3)
$$
as representations of $SL(2,\BB C)$ (see, for example, \cite{GGV}).
Then the map $\gamma$ switches the discrete and continuous components:
$$
L^2_{cont} \bigl( SU(1,1) \bigr) \simeq L^2_{discr}(\widetilde{H}^3),
\qquad
L^2_{discr} \bigl( SU(1,1) \bigr) \simeq L^2_{cont}(\widetilde{H}^3).
$$
This explains the purpose of the Cayley transform -- it is easier to find the
projector onto the discrete component than onto the continuous one!
Once the projector onto $L^2_{discr}(\widetilde{H}^3)$ is found
we can pull it back to $\HR$ to get $\operatorname{P}^{cont}$.
We expect it to have the form
$$
(\operatorname{P}^{cont}) f(X_0) = \frac 1 {2\pi^2} \int_{C_{cont}}
\frac{(Z-X_0)^{-1}}{N(Z-X_0)} \cdot Dz \cdot f(Z),
$$
where $C_{cont}$ is a certain three cycle in $\Omega$.

The integral formulas for $\operatorname{P}^{discr}$ and
$\operatorname{P}^{cont}$ strongly suggest that the separation of the series
is a topological phenomenon!
Thus there is some underlying homology theory still waiting to be developed.
We hope that this geometric analytic realization of the separation
of the series will extend to higher rank groups.
Such a picture fits well into the Gelfand-Gindikin program initiated in
\cite{GG} which realizes representations of reductive groups $G$ in function
spaces of open domains in $G_{\BB C}$ (complexification of $G$).
This geometric relationship between quaternionic analysis and the separation
of the series for $SL(2, \BB R)$ will be the subject of our upcoming paper
\cite{FL2}.

\separate

\separate

\noindent
{\em E-mail address:} {\tt mlibine@indiana.edu}

\noindent
{\em Department of Mathematics, Indiana University,
Rawles Hall, 831 East 3rd St, Bloomington, IN 47405}

\end{document}